\documentclass[12pt]{article}
\usepackage{latexsym,amssymb,amsmath,geometry,cite,enumerate,float}
\newtheorem{theorem}{Theorem}
\newtheorem{lemma}[theorem]{Lemma}

\newtheorem{corollary}[theorem]{Corollary}

\usepackage{lineno}
\usepackage{setspace}

\begin{document}

%\linenumbers
\onehalfspace

\title{Constant Threshold Intersection Graphs\\ of Orthodox Paths in Trees}

\author{Claudson Ferreira Bornstein$^1$
\and
Jos\'{e} Wilson Coura Pinto$^2$
\and
Dieter Rautenbach$^3$
\and
Jayme Luiz Szwarcfiter$^2$
}

\date{}

\maketitle

\begin{center}
$^1$ DCC-IM, Federal University of Rio de Janeiro, 
Rio de Janeiro, Brazil\\ 
\texttt{cfb@cos.ufrj.br}\\[3mm]
$^2$ PESC-COPPE, Federal University of Rio de Janeiro,
Rio de Janeiro, Brazil\\
\texttt{jwcoura,jayme@cos.ufrj.br}\\[3mm]
$^3$ Institute of Optimization and Operations Research,
Ulm University, Germany\\
\texttt{dieter.rautenbach@uni-ulm.de}
\end{center}

\begin{abstract}
A graph $G$ belongs to the class ${\rm ORTH}[h,s,t]$ for integers $h$, $s$, and $t$
if there is a pair $(T,{\cal S})$,
where $T$ is a tree of maximum degree at most $h$,
and ${\cal S}$ is a collection $(S_u)_{u\in V(G)}$ of subtrees $S_u$ of maximum degree at most $s$ of $T$, one for each vertex $u$ of $G$,
such that,
for every vertex $u$ of $G$, all leaves of $S_u$ are also leaves of $T$, and,
for every two distinct vertices $u$ and $v$ of $G$, 
the following three properties are equivalent:
\begin{enumerate}[(i)]
\item $u$ and $v$ are adjacent.
\item $S_u$ and $S_v$ have at least $t$ vertices in common.
\item $S_u$ and $S_v$ share a leaf of $T$.
\end{enumerate}
The class ${\rm ORTH}[h,s,t]$ was introduced by Jamison and Mulder.

Here we focus on the case $s=2$, which is closely related to the well-known VPT and EPT graphs.
We collect general properties of the graphs in ${\rm ORTH}[h,2,t]$,
and provide a characterization in terms of tree layouts.
Answering a question posed by Golumbic, Lipshteyn, and Stern,
we show that ${\rm ORTH}[h+1,2,t]\setminus {\rm ORTH}[h,2,t]$ is non-empty for every $h\geq 3$ and $t\geq 3$.
We derive decomposition properties, 
which lead to efficient recognition algorithms for the graphs in ${\rm ORTH}[h,2,2]$ for every $h\geq 3$.
Finally, we give a complete description of the graphs in ${\rm ORTH}[3,2,2]$, 
and show that the graphs in ${\rm ORTH}[3,2,3]$ are line graphs of planar graphs.
\end{abstract}

{\small 

\begin{tabular}{lp{13cm}}
{\bf Keywords:} & Intersection graph; $(h,s,t)$-representation; orthodox $(h,s,t)$-representation; line graph; chordal graph
\end{tabular}
}

\section{Introduction}

Intersection graphs are a well studied topic \cite{mcmc,gotr}
and the intersection graphs of paths in trees have received special attention. 
In the present paper 
we study so-called orthodox representations 
with bounds on the maximum degree of the host tree 
as well as 
on the size of the intersections corresponding to adjacencies.
Such representations were introduced by Jamison and Mulder \cite{jamu1,jamu2}.

Before we give precise definitions and explain our own as well as related results,
we collect some standard notation and terminology.
We consider finite, undirected, and simple graphs as well as finite and undirected multigraphs, which are allowed to contain parallel edges and loops.
A {\it clique} in $G$ is a complete subgraph of $G$.
For a tree $T$, let ${\cal L}(T)$ be the set of {\it leaves} of $T$,
which are the vertices of $T$ of degree at most $1$.
Let $L(H)$ be the {\it line graph} of some multigraph $H$,
whose vertex set $V(L(H))$ is the edge set $E(H)$ of $H$,
and in which two distinct vertices $u$ and $v$ of $L(H)$ are adjacent 
if and only if the edges $u$ and $v$ of $H$ intersect.
Two distinct vertices $u$ and $v$ of a graph $G$ are {\it twins} in $G$ if $N_G[u]=N_G[v]$, 
and, if $G$ has no twins, then it is {\it twin-free}.

The following notions were formalized by Jamison and Mulder \cite{jamu1,jamu2,jamu3}.
For positive integers $h$, $s$, and $t$,
an {\it $(h,s,t)$-representation} of a graph $G$ is a pair $(T,{\cal S})$,
where $T$ is a tree of maximum degree at most $h$,
and ${\cal S}$ is a collection $(S_u)_{u\in V(G)}$ of subtrees $S_u$ of maximum degree at most $s$ of $T$, one for each vertex $u$ of $G$,
such that two distinct vertices $u$ and $v$ of $G$ are adjacent 
if and only if $S_u$ and $S_v$ have at least $t$ vertices in common.
An $(h,s,t)$-representation $(T,{\cal S})$ of $G$ with ${\cal S}=(S_u)_{u\in V(G)}$ is {\it orthodox} if, 
for every vertex $u$ of $G$, all leaves of $S_u$ are also leaves of $T$, and,
for every two distinct vertices $u$ and $v$ of $G$, 
the following three properties are equivalent:
\begin{enumerate}[(i)]
\item $u$ and $v$ are adjacent.
\item $S_u$ and $S_v$ have at least $t$ vertices in common.
\item $S_u$ and $S_v$ share a leaf of $T$.
\end{enumerate}
Let $[h,s,t]$ and ${\rm ORTH}[h,s,t]$ be the classes of graphs that have an $(h,s,t)$-representation
and an orthodox $(h,s,t)$-representation, respectively.
If no upper bound on the maximum degree of the host $T$ is imposed, 
we replace $h$ with $\infty$.
Similarly, if no upper bound on the maximum degree of the subtrees in ${\cal S}$ is imposed, 
we replace $s$ with $\infty$.
Note that the classes $[h,s,t]$ and ${\rm ORTH}[h,s,t]$ are hereditary,
that is, closed under taking induced subgraphs.
By iteratively removing irrelevant leaves of the host tree $T$ of some orthodox $(h,s,t)$-representation $(T,{\cal S})$,
one may assume that every leaf of $T$ is also a leaf of some tree in ${\cal S}$.

Using this terminology, 
Gavril's famous result \cite{ga1} states that the class of chordal graphs coincides with $[\infty,\infty,1]$.
Jamison and Mulder \cite{jamu2,jamu3} attribute to McMorris and Scheinerman \cite{mcsc}
the insight that $[\infty,\infty,1]={\rm ORTH}[3,3,1]={\rm ORTH}[3,3,2]$.
In \cite{jamu2} they collect several properties of $[3,3,3]$ and ${\rm ORTH}[3,3,3]$.
The well studied 
vertex and edge intersection graphs of paths in trees \cite{ga2,goja1,goja2}, 
also known as VPT-graphs and EPT-graphs,
coincide with $[\infty,2,1]$ and $[\infty,2,2]$, respectively.
Golumbic and Jamison \cite{goja1} have shown that deciding whether a given graph belongs to $[3,2,1]$ is NP-complete.
Alc\'on, Gutierrez, and Mazzoleni \cite{alguma1} strengthened this result and generalized it for every $h\geq 3$.
In \cite{alguma2} they study the forbidden induced subgraphs of $[h,2,1]$.
Golumbic, Lipshteyn, and Stern \cite{golist} 
study the classes $[h,2,t]$ and ${\rm ORTH}[h,2,t]$ in detail.
In particular, they show that 
${\rm ORTH}[\infty,2,1]={\rm ORTH}[3,2,1]={\rm ORTH}[3,2,2]$,
and that ${\rm ORTH}[\infty,2,1]$ is a proper subclass of ${\rm ORTH}[\infty,2,2]$.
Furthermore, they ask whether ${\rm ORTH}[\infty,2,t]$ and ${\rm ORTH}[3,2,t]$ coincide.

\medskip

\noindent In the present paper we study the classes ${\rm ORTH}[h,2,t]$.

For $h\leq 2$, these classes are rather simple.
In fact, for every graph $G$ in ${\rm ORTH}[2,2,1]$, 
the vertex set $V(G)$ of $G$ can be partitioned into three cliques $A$, $B$, and $C$,
such that $G$ contains all edges between $A$ and $B$, all edges between $B$ and $C$, but no edge between $A$ and $C$,
that is, the only connected twin-free graph in ${\rm ORTH}[2,2,1]$ is $P_3$.
Furthermore, if $t\geq 2$, then the graphs in ${\rm ORTH}[2,2,t]$ consist of one clique and some isolated vertices.
Hence, the smallest interesting value for $h$ is $3$.

In the second section we collect some general properties of the graphs in ${\rm ORTH}[h,2,t]$.
Our main result is a characterization in terms of tree layouts whose precise definition will be given later.
Using this characterization, we are able to answer the above-mentioned question by Golumbic, Lipshteyn, and Stern.
We also derive some decomposition properties, 
which lead to efficient recognition algorithms for the graphs in ${\rm ORTH}[h,2,t]$ for every $h\geq 3$ and $t\in \{ 1,2\}$;
contrasting the above hardness results.
In the third section we consider the classes ${\rm ORTH}[3,2,2]$ and ${\rm ORTH}[3,2,3]$ in more detail,
and give a complete structural description of the first one.
We conclude with some open problems motivated by our research.

\section{General properties of ${\rm ORTH}[h,2,t]$}

In this section we collect more general properties of the classes ${\rm ORTH}[h,2,t]$,
and derive important structural consequences.
Our first result closely ties these classes to line graphs.

\begin{theorem}\label{theorem1}
Let $(T,{\cal S})$ be an orthodox $(h,2,t)$-representation of a graph $G$ with $h\geq 3$ and $t\geq 1$
such that, for every leaf $x$ of $T$, 
there is some vertex $u$ of $G$
with $x\in V(S_u)$.

The graph $G$ is the line graph of a multigraph $H$ without loops,
and, if $G$ is twin-free, then $H$ is a graph. 
Furthermore, if $G$ is a connected twin-free graph of order at least $4$,
and $H$ has no isolated vertices, then 
\begin{itemize}
\item $H$ is unique up to isomorphism,
\item there is a bijection $\phi:V(H)\to {\cal L}(T)$, and
\item two distinct vertices $x$ and $y$ of $H$ are adjacent in $H$
if and only if 
${\cal S}$ contains the path in $T$ between $\phi(x)$ and $\phi(y)$.
\end{itemize}
\end{theorem}
{\it Proof:} Let ${\cal S}=(S_u)_{u\in V(G)}$.
By definition, 
for every vertex $u$ of $G$,
the subtree $S_u$ is a path between leaves of $T$.

Suppose that, for some vertex $u$ of $G$, 
the path $S_u$ consists of only a single leaf, say $x$, of $T$.
Let $y$ be the neighbor of $x$ in $T$.
Let $T'$ arise from $T$ by adding the two new vertices $x'$ and $x''$
as well as the two new edges $xx'$ and $xx''$.
Now,
\begin{itemize}
\item replacing $S_v$ by the path $x'xx''$ for every vertex $v$ of $G$ for which $S_v$ consists only of $x$, and,
\item extending $S_v$ by the edge $xx'$ for every vertex $v$ of $G$ for which $S_v$ contains the edge $xy$ of $T$
\end{itemize}
yields an alternative orthodox $(h,2,t)$-representation of $G$ using $T'$ as host tree.
Possibly applying this transformation several times,
we may assume that every path in ${\cal S}$ has positive length.

By definition, 
for every leaf $x$ of $T$, 
the set 
$C_x=\{ u\in V(G):x\in V(S_u)\}$
is a clique in $G$.
Since, for every vertex $u$ of $G$, 
the subtree $S_u$ is a path between two distinct leaves of $T$,
every vertex of $G$ belongs to exactly two of the cliques in the collection $(C_x)_{x\in {\cal L}(T)}$.
Furthermore, for every edge $uv$ of $G$,
the two subtrees $S_u$ and $S_v$ share a leaf, say $x$, of $T$,
which implies that $u$ and $v$ both belong to $C_x$.
By results of Krausz \cite{kr} and of Bermond and Meyer \cite{beme},
this implies that $G$ is the line graph of some multigraph $H$ without loops.
Since parallel edges in $H$ correspond to twins in $G$,
if $G$ is twin-free, then $H$ is a graph.

Now, let $G$ be twin-free, connected, and of order at least $4$.
By a result of Whitney \cite{wh}, the graph $H$ is uniquely determined.
Let $H'$ be the graph with $V(H')={\cal L}(T)$
in which two distinct vertices $x$ and $y$ are adjacent
if and only if 
${\cal S}$ contains the path in $T$ between $x$ and $y$.
Since every leaf of $T$ is also a leaf of some $S_u$, 
the graph $H'$ has no isolated vertex.
By definition, the graph $G$ is isomorphic to $L(H')$,
and, by Whitney's result, the graphs $H'$ and $H$ are isomorphic.
$\Box$

\medskip

\noindent If the graph $G'$ arises from a graph $G$ 
by identifying all pairs of twins in $G$,
then $G'$ is twin-free.
Furthermore, it follows easily from the definition that 
$G$ is in ${\rm ORTH}[h,2,t]$
if and only if $G'$ is.

For $i\in [2]$, let $(T_i,{\cal S}_i)$ be an orthodox $(h,2,t)$-representation of a graph $G_i$ for some $h\geq 3$ and $t\geq 1$,
where $G_1$ and $G_2$ as well as $T_1$ and $T_2$ are disjoint.
Let $T$ be the tree that arises from the disjoint union of $T_1$ and $T_2$
by subdividing one edge of $T_i$ containing a leaf of $T_i$ with a new vertex $t_i$ for $i\in [2]$,
and adding the edge $t_1t_2$.
Since $h\geq 3$, applying the same subdivisions to the trees in ${\cal S}_1\cup {\cal S}_2$,
it follows that there is an orthodox $(h,2,t)$-representation of the disjoint union of $G_1$ and $G_2$
using $T$ as host tree.

In view of these observations and Theorem \ref{theorem1},
in order to understand the classes ${\rm ORTH}[h,2,t]$,
it suffices to consider connected twin-free line graphs $G$ 
of order at least $4$ as well as connected graphs $H$ with $L(H)=G$.

\medskip

\noindent The next result is our central characterization of the graphs in ${\rm ORTH}[h,2,t]$.

\begin{theorem}\label{theorem2}
Let $G$ be a connected twin-free line graph of order at least $4$, 
and let $H$ be a connected graph with $L(H)=G$. 

The graph $G$ is in ${\rm ORTH}[h,2,t]$
for some $h\geq 3$ and $t\geq 1$
if and only if there is a tree $T$ 
whose internal vertices all have degree at most $h$
such that $V(H)={\cal L}(T)$, and,
for every two independent edges $xy$ and $x'y'$ of $H$,
the two paths in $T$
between $x$ and $y$ and between $x'$ and $y'$
share at most $t-1$ vertices.
\end{theorem}
{\it Proof:} Let $G$ be in ${\rm ORTH}[h,2,t]$.
As shown in the proof of Theorem \ref{theorem1},
there is an orthodox $(h,2,t)$-representation $(T,{\cal S})$ of $G$
such that every path in ${\cal S}=(S_u)_{u\in V(G)}$ has positive length.
By definition, all internal vertices of $T$ have degree at most $h$.
By Theorem \ref{theorem1}, we may assume that $V(H)={\cal L}(T)$,
and that two distinct vertices $x$ and $y$ of $H$ are adjacent in $H$
if and only if 
${\cal S}$ contains the path in $T$ between $x$ and $y$.
If $xy$ and $x'y'$ are two independent edges of $H$,
then there are two vertices $u$ and $v$ in $G$ such that 
$S_u$ is the path in $T$ between $x$ and $y$ 
and
$S_v$ is the path in $T$ between $x'$ and $y'$.
Since $S_u$ and $S_v$ share no leaf, 
the two vertices $u$ and $v$ are not adjacent in $G$,
which implies that $S_u$ and $S_v$ share at most $t-1$ vertices.

Now, let $T$ be as in the statement.
Let $T'$ arise from $T$ by subdividing each edge incident with a leaf of $T$ exactly $t-2$ times.
Note that $T'$ still has maximum degree at most $h$,
and that,
for every two independent edges $xy$ and $x'y'$ of $H$,
the two paths in $T'$
between $x$ and $y$ and between $x'$ and $y'$
share at most $t-1$ vertices.
Let ${\cal S}=(S_{xy})_{xy\in E(H)}$,
where $S_{xy}$ is the path in $T'$ between the leaves $x$ and $y$ of $T'$.
Let $u$ and $v$ be two distinct vertices of $G$.
Let $u=xy$ and $v=x'y'$, 
where $xy$ and $x'y'$ are the corresponding edges of $H$.
Now, $u$ and $v$ are adjacent in $G$
if and only if the edges $xy$ and $x'y'$ are not independent 
if and only if $S_{xy}$ and $S_{x'y'}$ share a leaf of $T'$
if and only if $S_{xy}$ and $S_{x'y'}$ have at least $t$ vertices in common.
This implies that $(T',{\cal S})$ is an orthodox $(h,2,t)$-representation of $G$.
$\Box$

\medskip

\noindent The relation between $H$ and $T$ in the previous theorem is crucial for all our considerations, 
and we introduce some corresponding terminology.

If $H$ is a graph and $T$ is a tree such that 
\begin{itemize}
\item the maximum degree of $T$ is at most $h$,
\item $V(H)={\cal L}(T)$, and,
\item for every two independent edges $xy$ and $x'y'$ of $H$,
the two paths in $T$
between $x$ and $y$ and between $x'$ and $y'$
share at most $t-1$ vertices
\end{itemize}
for some integers $h\geq 3$ and $t\geq 1$,
then $T$ is an {\it $(h,t)$-tree layout} of $H$.

As we see now, vertices of degree $2$ are not essential within tree layouts.

\begin{corollary}\label{corollary1}
Let $G$ be a connected twin-free line graph of order at least $4$, 
and let $H$ be a connected graph with $L(H)=G$. 

The graph $G$ is in ${\rm ORTH}[3,2,t]$ for some $t\geq 1$
if and only if there is a tree $T$ 
whose internal vertices all have degree exactly $3$
such that $V(H)={\cal L}(T)$, and,
for every two independent edges $xy$ and $x'y'$ of $H$,
the two paths in $T$
between $x$ and $y$ and between $x'$ and $y'$
share at most $t-1$ vertices.
\end{corollary}
{\it Proof:} In view of Theorem \ref{theorem2},
it suffices to argue that no internal vertex of $T$ needs degree $2$.
Therefore, let $T$ be an $(h,t)$-tree layout of $H$.
If some vertex $b$ of $T$ 
has exactly two neighbors $a$ and $c$ within $T$,
then it is easy to see that $T'=T-b+ac$ is still an $(h,t)$-tree layout of $H$.
Iterating this transformation, it is possible to eliminate all internal vertices of $T$ that are of degree $2$. $\Box$

\medskip

\noindent Note that every tree with $n$ leaves whose internal vertices all have degree $3$, has order exactly $2n-2$.
Therefore, if $G$, $H$, and $T$ are as in the statement of Corollary \ref{corollary1}, then subdividing the edges incident with leaves of $T$ at most $t-2$ times, and defining ${\cal S}$ as in the proof of Theorem \ref{theorem2} yields an orthodox $(3,2,t)$-representation of $G$
whose underlying tree has order between $2n(H)-2$ and $tn(H)-2$.
Hence, by Theorem \ref{theorem2}, the minimum order of an underlying tree in any orthodox $(3,2,t)$-representation of $G$ lies between these two bounds. 

Our next goal is to answer the question posed by Golumbic, Lipshteyn, and Stern \cite{golist}.
In fact, we show that line graphs of complete graphs of suitable orders distinguish the classes ${\rm ORTH}[h,2,t]$ for different values of $h$.
The next lemma is a simple exercise, and we include the proof for completeness.

\begin{lemma}\label{lemma3}
Let $h$ and $t$ be integers with $h\geq 3$ and $t\geq 3$.

If $T$ is a tree of maximum degree at most $h$ such that every two leaves of $T$ have distance at most $t$, then 
$$
|{\cal L}(T)|\leq 
\begin{cases}
2(h-1)^{\left(\frac{t-1}{2}\right)} & \mbox{, if $t$ is odd, and}\\
h(h-1)^{\left(\frac{t}{2}-1\right)} & \mbox{, if $t$ is even.}
\end{cases}
$$
Furthermore, these bounds are tight for all considered values of $h$ and $t$.
\end{lemma}
{\it Proof:} In view of the desired bound, we may assume that $T$ has two leaves at distance $t$.
Let $x(0)\ldots x(t)$ be a path between two such leaves.

First, let $t$ be odd. 
Rooting the two components of $T-x\left(\frac{t-1}{2}\right)x\left(\frac{t-1}{2}+1\right)$ 
in the two vertices $x\left(\frac{t-1}{2}\right)$ and $x\left(\frac{t-1}{2}+1\right)$
yields rooted $(h-1)$-ary trees $T_1$ and $T_2$ of depth $\frac{t-1}{2}$.
Clearly,
$$|{\cal L}(T)|=|{\cal L}(T_1)|+|{\cal L}(T_2)|\leq (h-1)^{\left(\frac{t-1}{2}\right)}+(h-1)^{\left(\frac{t-1}{2}\right)}$$
with equality if and only if 
$T_1$ and $T_2$ are full $(h-1)$-ary trees 
of depth $\frac{t-1}{2}$.

Next, let $t$ be even. 
Rooting the two components of $T-x\left(\frac{t}{2}\right)x\left(\frac{t}{2}+1\right)$ 
in the two vertices $x\left(\frac{t}{2}\right)$ and $x\left(\frac{t}{2}+1\right)$
yields rooted $(h-1)$-ary trees $T_1$ of depth $\frac{t}{2}$ and $T_2$ of depth $\frac{t}{2}-1$.
Clearly,
$$|{\cal L}(T)|=|{\cal L}(T_1)|+|{\cal L}(T_2)|\leq (h-1)^{\left(\frac{t}{2}\right)}+(h-1)^{\left(\frac{t}{2}-1\right)}=h(h-1)^{\left(\frac{t}{2}-1\right)}$$
with equality if and only if 
$T_1$ and $T_2$ are full $(h-1)$-ary trees 
of depths $\frac{t}{2}$ and $\frac{t}{2}-1$, respectively.
$\Box$ 

\medskip

\noindent The following result answers the question posed by Golumbic, Lipshteyn, and Stern \cite{golist} mentioned in the introduction.
For two integers $p$ and $q$, let $[p,q]$ denote the set of integers at least $p$ and at most $q$.

\begin{theorem}\label{theorem4}
Let $h$ and $t$ be integers with $h\geq 3$ and $t\geq 3$.

If $t$ is odd, then 
$$\Big\{ n\in \mathbb{N}: L(K_n)\in {\rm ORTH}[h+1,2,t]\setminus {\rm ORTH}[h,2,t]\Big\}
=\left[ 2(h-1)^{\left(\frac{t-1}{2}\right)}+1,2h^{\left(\frac{t-1}{2}\right)}\right],$$
and, if $t$ is even, then 
$$\Big\{ n\in \mathbb{N}: L(K_n)\in {\rm ORTH}[h+1,2,t]\setminus {\rm ORTH}[h,2,t]\Big\}
=\left[ h(h-1)^{\left(\frac{t}{2}-1\right)}+1,(h+1)h^{\left(\frac{t}{2}-1\right)}\right].$$
\end{theorem}
{\it Proof:} We only give details for odd $t$,
because the proof for even $t$ is analogous.
Since the considered classes are hereditary,
it suffices to show that
\begin{enumerate}[(i)]
\item $L\left(K_{2h^{\left(\frac{t-1}{2}\right)}}\right)\in {\rm ORTH}[h+1,2,t]$, and
\item 
$L\left(K_{2(h-1)^{\left(\frac{t-1}{2}\right)}}+1\right)\not\in {\rm ORTH}[h,2,t]$.
\end{enumerate}
By Lemma \ref{lemma3}, there is a tree of maximum degree at most $h+1$
with $2h^{\left(\frac{t-1}{2}\right)}$ leaves 
such that every two leaves have distance at most $t$.
By Theorem \ref{theorem2}, this implies (i).

Now, suppose that (ii) does not hold.
Again by Theorem \ref{theorem2}, 
there is a tree $T$ of maximum degree at most $h$
with $2h^{\left(\frac{t-1}{2}\right)}$ leaves such that,
for every four distinct leaves $u_1$, $v_1$, $u_2$, and $v_2$,
the two paths in $T$ between $u_1$ and $v_1$
and between $u_2$ and $v_2$ share at most $t-1$ vertices.
Let $u_1$ and $v_1$ be two leaves of $T$ with maximum distance $\ell$.
We assume that $T$ is chosen such that $\ell$ is as small as possible.
Let $u_1'$ and $v_1'$ be the two neighbors of $u_1$ and $v_1$, respectively.
Clearly, the vertices $u_1'$ and $v_1'$ are distinct.
By the choice of $T$,
the vertex $u_1'$ is adjacent to a leaf $u_2$ distinct from $u_1$,
and the vertex $v_1'$ is adjacent to a leaf $v_2$ distinct from $v_1$.
Considering the two independent edges $u_1v_1$ and $u_2v_2$ implies that $\ell\leq t$.
By Lemma \ref{lemma3}, this implies the contradiction
that $T$ has at most $2(h-1)^{\left(\frac{t-1}{2}\right)}$ leaves. $\Box$

\medskip

\noindent Next, we show that all subgraphs of the essential graphs $H$ with $L(H)\in {\rm ORTH}[h,2,t]$ 
have balanced separations of bounded order. 
 
\begin{theorem}\label{theorem5}
Let $G$ be a connected twin-free line graph of order at least $4$, 
and let $H$ be a connected graph with $L(H)=G$. 
Let $G$ be in ${\rm ORTH}[h,2,t]$ for some $h\geq 3$ and $t\geq 1$.

If $H'$ is a subgraph of $H$ of order at least $2$,
then there is a set $X$ of vertices of $H'$, and 
a partition of $V(H')$ into two sets $A$ and $B$
such that 
\begin{enumerate}[(i)]
\item $\frac{1}{h}n(H')\leq |A|,|B|\leq \frac{h-1}{h}n(H')$, 
\item $|X|\leq \max\left\{ 1,(h-1)^{(t-2)}\right\}$, and
\item $H'$ contains no edge between $A\setminus X$ and 
$B\setminus X$.
\end{enumerate}
Furthermore, given $H'$, the sets $X$, $A$, and $B$ can be found in polynomial time.
\end{theorem}
{\it Proof:} Let $T$ be as in Theorem \ref{theorem2},
that is, the tree $T$ is an $(h,t)$-tree layout of $H$.
Clearly, iteratively removing leaves of $T$ that are not vertices of $H'$
yields an $(h,t)$-tree layout $T'$ of $H'$.

If $t=1$, 
$u$ is some vertex of $T'$, 
$Y$ is the set of leaves of some components of $T'-u$,
and 
$Z$ is the set of leaves of the remaining components of $T'-u$,
then the properties of a $(h,1)$-tree layout 
imply the existence of a vertex $x$ in $Y\cup Z$ 
such that all edges of $H'$ between $Y$ and $Z$ 
are incident with $x$,
because otherwise $H'$ contains two independent edges between $Y$ and $Z$
such that the corresponding paths in $T'$ share $u$.
Similarly, if $t\geq 2$, $P:u_1\ldots u_t$ is a path of order $t$ in $T'$, 
$Y$ is the set of leaves of $T'$ that lie in the same component of $T'-E(P)$ as $u_1$, and
$Z$ is the set of leaves of $T'$ that lie in the same component of $T'-E(P)$ as $u_t$,
then there is a vertex $x$ in $Y\cup Z$ 
such that all edges of $H'$ between $Y$ and $Z$ 
are incident with $x$.

Let $r$ be any internal vertex of $T'$, and root $T'$ in $r$.
Let $a$ be a vertex of $T'$ of maximum distance from $r$
such that at least $\frac{1}{h}n(H')$ descendants of $a$ are leaves.
Since $T'$ has $n(H')$ leaves and $r$ has degree at most $h$,
the vertex $a$ is not $r$.
Since $a$ has at most $h-1$ children, 
at most $\frac{h-1}{h}n(H')$ descendants of $a$ are leaves.
Let $b$ be the parent of $a$.
Let $A$ be the set of leaves of $T'$ that lie in the same component of $T'-ab$ as $a$, and 
let $B$ be the set of leaves of $T'$ that lie in the same component of $T'-ab$ as $b$.
Clearly, (i) holds.

If $t=1$, then the above observation implies the existence of a single vertex $x$ such that all edges of $H'$ between $A$ and $B$ are incident with $x$, and, (ii) and (iii) follow.
Now, let $t\geq 2$.
Rooting the component $T'_b$ of $T'-ab$ that contains $b$ in the vertex $b$,
and considering all leaves of $T'_b$ at depth at most $t-2$
as well as all non-leaf vertices of $T'_b$ at depth exactly $t-2$,
it follows that $B$ can be partitioned into $k\leq (h-1)^{t-2}$ sets $B_1,\ldots,B_k$
such that, for every $i\in [k]$, 
if $B_i$ contains not only one vertex, 
then there is a path $P_i$ of order $t$ 
such that every path in $T'$ between a leaf in $A$
and a leaf in $B_i$ has $P_i$ as a subpath.
If $B_i$ contains only one vertex, 
then trivially all edges of $H'$ between $A$ and $B_i$ are incident with only one vertex in $A\cup B_i$.
If $B_i$ contains not only one vertex, 
then the above observation also implies that 
all edges of $H'$ between $A$ and $B_i$ are incident with only one vertex in $A\cup B_i$.
Altogether, it follows that there is a set $X$ of at most $k$ vertices of $H'$
such that all edges between $A$ and 
$B=B_1\cup \ldots \cup B_k$ 
are incident with a vertex in $X$, and, (ii) and (iii) follow.

It remains to explain, how to determine suitable sets $X$, $A$, and $B$ efficiently given $H'$.
Therefore, let $p=\max\left\{ 1,(h-1)^{(t-2)}\right\}$.
If $p\geq \frac{h-1}{h}n(H')$, 
then choosing $A$ as any set of $\left\lceil\frac{1}{h}n(H')\right\rceil$ 
vertices of $H'$, and $B=X$ as $V(H')\setminus A$
satisfies (i), (ii), and (iii).
Now, let $p<\frac{h-1}{h}n(H')$.
For any specific set $X$ of $p$ vertices of $H'$,
we explain how to decide whether sets $A$ and $B$ with (i), (ii), and (iii) exist.
Therefore, let $X$ be such a set.
Let $K_1,\ldots,K_\ell$ be the components of $H'-X$.
If some component $K_i$ with $i\in [\ell]$ has order more than $\frac{h-1}{h}n(H')$,
then the desired sets $A$ and $B$ do not exist.
If some component $K_i$ with $i\in [\ell]$ has order between $\frac{1}{h}n(H')$ and $\frac{h-1}{h}n(H')$,
then choosing $A=V(K_i)$ and $B=V(H')\setminus A$ has the desired properties. 
Finally, 
if all components $K_i$ with $i\in [\ell]$ have order less than $\frac{1}{h}n(H')$,
then forming a union of suitably many of their vertex sets 
yields a set $A$ of order between $\frac{1}{h}n(H')$ and $\frac{h-1}{h}n(H')$, and defining $B$ as above
leads to sets with the desired properties.
Altogether, considering the $O(n(H')^p)$ many choices for $X$,
suitable sets can be determined in polynomial time,
which completes the proof.
$\Box$

\medskip

\noindent The following immediate consequence of Theorem \ref{theorem5}
might be useful in order to devise efficient recognition algorithms for the graphs in ${\rm ORTH}[h,2,t]$.

\begin{corollary}\label{corollary2}
If $G$ is a connected twin-free line graph of order at least $4$ in ${\rm ORTH}[h,2,t]$ for some $h\geq 3$ and $t\geq 1$, and $H$ is a connected graph with $L(H)=G$, then the treewidth of $H$ is bounded as a function of $h$ and $t$.
\end{corollary}
{\it Proof:} This follows easily from Theorem \ref{theorem5} 
and a result of Dvo\v{r}\'ak and Norin \cite{dvno}. $\Box$

\medskip

\noindent Our final goal in this section are efficient recognition algorithms for the graphs in ${\rm ORTH}[h,2,1]$ and ${\rm ORTH}[h,2,2]$
based on Theorem \ref{theorem5}.

\begin{lemma}\label{lemma4}  
Let $H$ be a graph.
Let the sets $A$ and $B$ partition $V(H)$.
Let $a$ in $A$ and let $b$ in $B$ be such that 
every edge of $H$ between $A$ and $B$ is incident with $a$,
and $b$ is a neighbor of $a$.

The graph $H$ has an $(h,t)$-tree layout
for some $h\geq 3$ and $t\geq 1$
if and only if
the two graphs $H_A=H[A\cup \{ b\}]$ and $H_B=H[B\cup \{ a\}]$ have $(h,t)$-tree layouts.
\end{lemma}
{\it Proof:} If $H$ has an $(h,t)$-tree layout,
then so does every induced subgraph, which implies the necessity.
For the sufficiency, assume that 
$T_A$ and $T_B$ are $(h,t)$-tree layouts of $H_A$ and $H_B$,
respectively.
Note that $T_A$ and $T_B$ share exactly the two vertices $a$ and $b$.
If the tree $T$ arises from the disjoint union of $T_A$ and $T_B$,
where we distinguish the two copies of $a$ and $b$ within $T_A$ and $T_B$,
by adding an edge between the copy of $b$ in $T_A$ and the copy of $a$ in $T_B$,
then it follows easily that $T$ is an $(h,t)$-tree layout of $H$,
which completes the proof. $\Box$

\begin{corollary}\label{corollary3}
For every two integers $h\geq 3$ and $t\in \{ 1,2\}$, 
the graphs in ${\rm ORTH}[h,2,t]$ can be recognized in polynomial time.
\end{corollary}
{\it Proof:} Let $G$ be a given graph for which we want to decide whether is belongs to ${\rm ORTH}[h,2,t]$.
As observed after Theorem \ref{theorem1},
we may assume that $G$ is a connected twin-free line graph of order at least $4$.
Using for instance the algorithms in \cite{desi},
we can efficiently determine the unique connected graph $H$ 
with $L(H)=G$.
Clearly, $n(H)\leq m(H)+1=n(G)+1$.
By Theorem \ref{theorem2},
we need to decide whether $H$ has an $(h,t)$-tree layout.

Note that, for $t\in \{ 1,2\}$,
the set $X$ in Theorem \ref{theorem5} contains at most one vertex.
Furthermore, sets $X$, $A$, and $B$ with (i), (ii), and (iii) can be found efficiently for every subgraph $H'$ of $H$ of order at least $2$.
Note that the graphs $H_A$ and $H_B$ considered in Lemma \ref{lemma4} have orders $|A|+1$ and $|B|+1$, respectively.
Let $n_0$ be such that $\frac{h-1}{h}n+1\leq \frac{h}{h+1}n$ for $n\geq n_0$. 
Note that, if $H$ has order at least $n_0$, 
then $H_A$ and $H_B$ both have orders at most $\frac{h}{h+1}n(H)$.
Therefore, 
iteratively applying Theorem \ref{theorem5} yields 
$k\leq 2^{\lceil \log_{(h+1)/h}(n/n_0)\rceil}$ many graphs $H_1,\ldots,H_k$,
each of order at most $n_0$, 
such that 
$H$ has an $(h,t)$-tree layout 
if and only each $H_i$ an $(h,t)$-tree layout for every $i\in [k]$.
Clearly, testing this property for these polynomially many graphs of bounded order can be done in polynomial time.
$\Box$

\section{${\rm ORTH}[3,2,2]$ and ${\rm ORTH}[3,2,3]$}

This section collects more specific properties of the classes ${\rm ORTH}[3,2,2]$ and ${\rm ORTH}[3,2,3]$.
As shown by Golumbic, Lipshteyn, and Stern \cite{golist},
we have ${\rm ORTH}[\infty,2,1]={\rm ORTH}[3,2,1]={\rm ORTH}[3,2,2]$,
which implies ${\rm ORTH}[h,2,1]={\rm ORTH}[3,2,2]$ for every $h\geq 3$.

\begin{theorem}\label{theorem6}
Let $G$ be a connected twin-free line graph of order at least $4$, 
and let $H$ be a connected graph with $L(H)=G$. 

The graph $G$ is in ${\rm ORTH}[3,2,2]$
if and only if 
all blocks of $H$ are of order at most $3$.
\end{theorem}
{\it Proof:} 
In order to show the necessity, 
we assume that $G$ is in ${\rm ORTH}[3,2,1]={\rm ORTH}[3,2,2]$
but that some block of $H$ has order at least $4$.
This implies that $H$ has some cycle $C$ of length $\ell$ at least $4$ as a subgraph.
Hence, $G$ contains the induced cycle $L(C)$ of length $\ell$.
Nevertheless, by Gavril's result \cite{ga1}, the graphs in ${\rm ORTH}[3,2,1]$ are chordal,
which is a contradiction.

In order to show the sufficiency, 
we assume that all blocks of $H$ are of order at most $3$.
Since $K_3$ has a $(3,1)$-tree layout, it follows easily by an inductive argument similar to the proof of Lemma \ref{lemma4}
that $H$ has a $(3,1)$-tree layout.
Suppose, for instance, 
that $H$ arises from a smaller graph $H'$ 
containing a vertex $u$ 
by adding two vertices $v$ and $w$, 
and three new edges $uv$, $uw$, and $vw$, 
that is, $H$ arises from $H'$ by attaching one new $K_3$ block.
If $T'$ is a $(3,1)$-tree layout of $H'$,
then subdividing the edge of $T'$ incident with $u$
by a new vertex $x$,
adding three more vertices $y$, $v$, and $w$,
and adding three more edges $xy$, $vy$, and $wy$
yields a $(3,1)$-tree layout of $H$.
By Theorem \ref{theorem2}, this implies that $G$ is in ${\rm ORTH}[3,2,1]={\rm ORTH}[3,2,2]$.
$\Box$

\begin{lemma}\label{lemma1}
If $e$ and $f$ are two independent edges of $K_5$, $H$ is a subdivision of $K_5-\{ e,f\}$, and $G=L(H)$,
then $G$ is not in ${\rm ORTH}[3,2,3]$.
\end{lemma}
{\it Proof:} Let $H_0=K_5-\{ e,f\}$.
We denote the five vertices of $H_0$ by $u_1,\ldots,u_5$
without specifying which two edges are missing.
Note that $G$ is a connected twin-free line graph of order at least $4$.
For a contradiction, suppose that $G$ is in ${\rm ORTH}[3,2,3]$.
Let $T$ be as in Corollary \ref{corollary1} for $t=3$.
Let $r$ be any internal vertex of $T$, and root $T$ in $r$.
Let $s$ be a vertex of $T$ of maximum distance from $r$
such that at least two descendants, say $u_1$ and $u_2$, 
of $s$ within $T$ are vertices of $H_0$.
Since every internal vertex of $T$ has degree $3$,
exactly two descendants of $s$ are vertices of $H_0$,
which implies that $s$ is not $r$.
Let $t$ be the parent of $s$.
Let $s'$ and $s''$ be the two neighbors of $t$ distinct from $s$.
Let $S$, $S'$, and $S''$ be the vertex sets of the three components of $T-t$
that contain $s$, $s'$, and $s''$, respectively.

By the pigeonhole principle and by symmetry,
we may assume that $u_4$ and $u_5$ lie in $S'$, 
and that $u_1u_4$ and $u_2u_5$ are edges of $H_0$.
Note that every edge $uv$ of $H_0$ corresponds to a path $P(uv)$ in $H$
between $u$ and $v$ whose internal vertices are all of degree $2$,
and that the edges of this path correspond to leaf to leaf paths in $T$
such that, for independent edges, the corresponding paths share at most $2$ vertices.

Since the two paths $P(u_1u_4)$ and $P(u_2u_5)$ are between a vertex in $S$ and a vertex in $S'$, 
each contains 
\begin{itemize}
\item an edge between $S$ and $S'$, or
\item an edge between $S$ and $S''$ as well as an edge between $S'$ and $S''$.
\end{itemize}
Since any edge of $P(u_1u_4)$ is disjoint from any edge of $P(u_2u_5)$,
the structure of $T$ implies that we may assume that 
$P(u_1u_4)$ contains an edge $e_1$ between $S$ and $S'$,
and that 
$P(u_2u_5)$ contains an edge $e_2$ between $S$ and $S''$ 
as well as an edge $e_2'$ between $S'$ and $S''$.
If $u_1$ is adjacent to $u_5$ in $H_0$, then, in view of $e_2$, the path $P(u_1u_5)$ contains an edge between $S$ and $S'$,
and,
if $u_2$ is adjacent to $u_4$ in $H_0$, then, in view of $e'_2$, the path $P(u_2u_4)$ contains an edge between $S$ and $S'$.
Since every edge of $P(u_1u_5)$ is disjoint from every edge of $P(u_2u_4)$,
this implies that $u_1u_5$ or $u_2u_4$ is one of the two non-edges of $H_0$.
If $u_3$ is in $S''$, then, in view of $e_2$ and $e_2'$, 
the vertex $u_3$ can not be adjacent to $u_1$ or $u_4$ in $H_0$, 
which is a contradiction.
Hence, $u_3$ is in $S'$.
In view of $e_1$ and $e_2'$, the vertex $u_3$ is not adjacent to $u_2$ in $H_0$.
Together with our earlier observation,
this implies that the two missing edges of $H_0$ are exactly $u_1u_5$ and $u_2u_3$.
In view of $e_2$, the path $P(u_1u_3)$ contains an edge $e_3$ between $S$ and $S'$.
In view of $e_2'$, the path $P(u_2u_4)$ contains an edge $e_4$ between $S$ and $S'$.
Now, the two paths in $T$ between the endpoints of $e_3$
as well as 
between the endpoints of $e_4$ share three vertices $s$, $t$, and $s'$,
which is a contradiction.
$\Box$

\medskip

\noindent The previous lemma has a suitable generalization for larger values of $h$ than $3$.
A similar proof also shows that $L(K_{2,5})$ does not lie in ${\rm ORTH}[3,2,3]$.

\begin{lemma}\label{lemma2}
If $H$ is a subdivision of $K_{3,3}$ and $G=L(H)$,
then $G$ is not in ${\rm ORTH}[3,2,3]$.
\end{lemma}
{\it Proof:} Let the two partite sets of $K_{3,3}$ be $\{ u_1,u_2,u_3\}$ and $\{ u'_1,u'_2,u'_3\}$.
Note that $G$ is a connected twin-free line graph of order at least $4$.
For a contradiction, suppose that $G$ is in ${\rm ORTH}[3,2,3]$.
Let $T$ be as in Corollary \ref{corollary1} for $t=3$.
Let $r$ be any internal vertex of $T$, and root $T$ in $r$.
Let $s$ be a vertex of $T$ of maximum distance from $r$
such that at least two descendants of $s$ within $T$ 
are vertices of $K_{3,3}$.
Again, exactly two descendants of $s$ are vertices of $K_{3,3}$.
Let $t$, $s'$, $s''$, $S$, $S'$, and $S''$ be as in the proof of Lemma \ref{lemma1}.

We consider some cases concerning the two vertices of $K_{3,3}$ in $S$
as well as the distribution of the remaining four vertices $K_{3,3}$ within $S'$ and $S''$.

\medskip

\noindent {\bf Case 1} {\it $u_1,u_2\in S$.}

\medskip

\noindent By symmetry, we may assume that $u_1',u_2'\in S'$.
Arguing as in the proof of Lemma \ref{lemma1},
we may assume, by symmetry, 
that $P(u_1u_1')$ contains an edge $e_1$ between $S$ and $S''$
as well as an edge $e_1'$ between $S'$ and $S''$,
and that $P(u_2u_2')$ contain an edge $e_2$ between $S$ and $S'$.
In view of $e_1$, the path $P(u_2u_3')$ contains no edge between $S$ and $S''$.
In view of $e_1'$, the path $P(u_2u_3')$ contains no edge between $S'$ and $S''$.
This implies that $u_3'\in S'$.
In view of $e_2$, the path $P(u_1u_3')$ contains no edge between $S$ and $S'$.
In view of $e_1'$, the path $P(u_1u_3')$ contains no edge between $S'$ and $S''$.
This implies a contradiction.

In view of Case 1, we may assume that the two vertices of $K_{3,3}$ in $S$ 
belong to different partite sets, that is, by symmetry, $u_1,u_1'\in S$.

\medskip

\noindent {\bf Case 2} {\it $u_1,u_1'\in S$, $u_2,u_2'\in S'$, and $u_3,u_3'\in S''$.}

\medskip

\noindent By symmetry, we may assume that 
$P(u_1u_2')$ contains an edge $e_1$ between $S$ and $S''$ as well as an edge $e_1'$ between $S'$ and $S''$, and that 
$P(u'_1u_2)$ contains an edge $e_2$ between $S$ and $S'$.
In view of $e_1$, the path $P(u_2u_3')$ contains no edge between $S$ and $S''$.
In view of $e_1'$, the path $P(u_2u_3')$ contains no edge between $S'$ and $S''$.
This implies a contradiction.

\medskip

\noindent {\bf Case 3} {\it $u_1,u_1'\in S$ and $u_2,u_2',u_3,u_3'\in S'$.}

\medskip

\noindent By symmetry, we may assume that 
$P(u_1u_2')$ contains an edge $e_1$ between $S$ and $S''$ as well as an edge $e_1'$ between $S'$ and $S''$, and that 
$P(u'_1u_2)$ contains an edge $e_2$ between $S$ and $S'$.
In view of $e_2$, the path $P(u_1u_3')$ contains no edge between $S$ and $S'$.
In view of $e_1'$, the path $P(u_1u_3')$ contains no edge between $S'$ and $S''$.
This implies a contradiction.

In view of Cases 2 and 3, 
we may assume, by symmetry, that $u_2\in S'$ and $u_2'\in S''$.
If $u_3\in S'$ and $u_3'\in S''$, then we can argue as in Case 1.
If $u_3'\in S'$ and $u_3\in S''$, then we can argue as in Case 2.
Hence, by symmetry, it suffices to consider the following final case.

\medskip

\noindent {\bf Case 4} {\it $u_1,u_1'\in S$, $u_2,u_3,u_3'\in S'$, and $u_2'\in S''$.}

\medskip

\noindent By symmetry, we may assume that 
$P(u_1u_3')$ contains an edge $e_1$ between $S$ and $S''$ as well as an edge $e_1'$ between $S'$ and $S''$, and that 
$P(u_1'u_3)$ contains an edge $e_2$ between $S$ and $S'$.
In view of $e_2$, the path $P(u_2u_2')$ contains no edge between $S$ and $S'$.
In view of $e_1'$, the path $P(u_2u_2')$ contains no edge between $S'$ and $S''$.
This implies a contradiction, and completes the proof. $\Box$

\begin{theorem}\label{theorem3}
If a connected twin-free graph $G$ of order at least $4$ is in ${\rm ORTH}[3,2,3]$, then $G$ is the line graph of a planar graph.
\end{theorem}
{\it Proof:} By Theorem \ref{theorem1},
there is a unique connected graph $H'$ with $G=L(H')$.
For a contradiction, suppose that $H'$ is not planar.
By a result of Kuratowski \cite{ku},
the graph $H'$ has a subgraph $H$ 
that is a subdivision of a graph $H_0$
such that $H_0$ is either $K_5$ or $K_{3,3}$.
Since $L(H)$ is an induced subgraph of $L(H')$,
the graph $L(H)$ is in ${\rm ORTH}[3,2,3]$.
Now, if $H_0$ is $K_5$ or $K_{3,3}$, 
then Lemma \ref{lemma1} or Lemma \ref{lemma2} 
implies the contradiction that $L(H)$ is not in ${\rm ORTH}[3,2,3]$,
respectively.
$\Box$

\medskip

\noindent Since Lemma \ref{lemma1} actually concerns subdivisions of $K_5-e$, the containment in Theorem \ref{theorem3} is proper.

\section{Conclusion}

The most natural open problems concern the structure 
of the graphs in the classes ${\rm ORTH}[3,2,3]$ and ${\rm ORTH}[h,2,2]$ for $h\geq 4$.
For ${\rm ORTH}[3,2,3]$, the complexity of the recognition is unknown.
In view of Corollary \ref{corollary2}, efficient recognition algorithms for all classes ${\rm ORTH}[h,2,t]$ seem possible.
Our results should have further algorithmic consequences.
If, for example, $G$ and $H$ are as in Theorem \ref{theorem2}, 
then the chromatic number of $G$ is either the maximum degree of $H$ or one more,
and, by Corollary \ref{corollary2}, these two cases can be distinguished efficiently.

\end{document}